\documentclass[11pt, final]{amsart}
\usepackage{rmwe}
\usepackage{geometry}
\usepackage[parfill=0pt]{parskip}    
\usepackage{graphicx}
\usepackage{amssymb, centernot}
\usepackage{epstopdf}
\usepackage{enumitem}
\setlist[enumerate,1]{ref=\emph{(\roman*)}}
\setlist[enumerate,1]{label=\emph{(\roman*)}}

\usetikzlibrary{decorations.pathmorphing}

\mynewtheorem{propdefn}{Proposition/Definition}{Prop/Defn}

\newcommand{\dashto}{\dashrightarrow}
\newcommand{\E}{\mathcal{E}}
\newcommand{\nin}{\notin}
\newcommand{\emp}{\emptyset}
\newcommand{\sgraph}{\Sigma}
\newcommand{\nice}{untangled}
\newcommand{\Nice}{Untangled}

\newcommand{\lo}{o}
\newcommand{\minimalalgname}{Minimal Stabilisation Algorithm}
\newcommand{\ev}[1]{#1_{\text{ev}}}

\makeatletter
\def\blfootnote{\xdef\@thefnmark{}\@footnotetext}
\makeatother

\DeclareMathOperator{\Pic}{Pic}
\DeclareMathOperator{\NS}{NS}
\DeclareMathOperator{\rk}{rk}
\DeclareMathOperator{\comp}{\mathfrak e}

\title[On the Stabilisation of Rational Surface Maps]{On the Stabilisation of Rational Surface Maps}
\author{\mylongname}
\address{Department of Mathematics, Statistics, and Computer Science, University of Illinois at Chicago, 851 S. Morgan Street, Chicago, IL 60607-7045}
\email{rbirkett@uic.edu}

\begin{document}

\begin{abstract}

The dynamics of a rational surface map $f : X \dashrightarrow X$ are easier to analyse when $f$ is `algebraically stable'.  Here we investigate when and how this condition can be achieved by conjugating $f$ with a birational change of coordinates.  We show that if this can be done with a birational morphism, then there is a minimal such conjugacy. For birational $f$ we also show that repeatedly lifting $f$ to its graph gives a stable conjugacy. Finally, we give an example in which $f$ can be birationally conjugated to a stable map, but the conjugacy cannot be achieved solely by blowing up.
 
 La dynamique d’une application rationnelle $f : X \dashrightarrow X$ sur une surface est plus simple \`a analyser lorsque $f$ est `alg\'ebriquement stable'. Dans cet article nous \'etudions comment la stabilit\'e peut \^etre r\'ealis\'ee en conjuguant $f$ par un changement de variable birationnel. Nous montrons que si cela peut \^etre r\'ealis\'ee avec une morphisme birationnelle, il existe une telle conjugaison minimale. Pour $f$ birationnelle, nous montrons aussi que l’on obtient une conjugaison stable par rel\`evement successif au graphe. Nous donnons enfin un exemple dans lequel $f$ peut \^etre conjugu\'ee birationnellement \`a une application stable, mais la conjugu\'ee ne peut pas \^etre obtenue uniquement par \'eclatement.

\end{abstract}

\maketitle
 \blfootnote{To appear in the Annales de la Facult\'e des Sciences de Toulouse.}

\vspace{-5mm}
\section{Introduction}

Let $f: X \dashto X$ be a rational map on a smooth projective surface over an algebraically closed field. Studying $f$ as a dynamical system is complicated by the fact that $f$ need not be continuously defined on all points of $X$. For example a rational map $f$ induces a natural pullback operator on curves $f^* : \Pic(X) \to \Pic(X)$, but this operator may not iterate well. We say $f$ is \emph{algebraically stable} iff $\forall n \in \N\ (f^*)^n = (f^n)^*$. This is equivalent to the geometric condition that $f$ has no \emph{destabilising orbits} \cite{FS} \cite[1.14]{DF}; i.e. an orbit of (closed) points $p, f(p)\dots, f^{n-1}(p)$ in $X$, for which $f^{-1}(p)$ and $f(f^{n-1}(p))$ are curves. It is natural to hope that blowing up the points in such an orbit will improve the situation. The main theme of this paper is to discuss the extent to which this actually works.

Assume for the rest of the paper that all surfaces are projective varieties over an algebraically closed field, and unless explicitly stated, also smooth. We write $\phi : X \dashto Y$ to indicate that $f$ is a rational map between surfaces, and we use a solid arrow $f : X \to Y$ when $f$ is a morphism.

\begin{defn}
 We write $\phi : (g, Y) \dashto (f, X)$ to indicate that $\phi : Y \dashto X$ is a birational map conjugating $f : X \dashto X$ to $g = \phi^{-1} \circ f \circ \phi : Y \dashto Y$. When $g : Y \dashto Y$ is algebraically stable, we say that $\phi$ \emph{stabilises} $f$.
\end{defn}

Diller-Favre \cite{DF} proved that for a birational map $(f, X)$ there is always a birational morphism $\pi : (g, Y) \to (f, X)$ which stabilises $f$. Not all (non-invertible) rational maps can be stabilised, however. Favre \cite{Fav} showed for example that many monomial maps on $\P^2$ cannot be stabilised by any birational conjugacy. See also \cite{HP}; for discussion on monomial maps in higher dimensions see \cite{JW, Lin}, and for the local case at normal surface singularities see \cite{Fav2, GR2}.

Algebraic stability is valuable in particular because the first dynamical degree becomes easy to compute. This can be defined by \[\lambda_1(f) = \lim_{n\to \infty} \|(f^n)^*\|^{1/n}.\]
For an $f : X \dashto X$ which is not algebraically stable, this computation is usually intractable. However, if one has a stabilisation $\phi : (g, Y) \dashto (f, X)$ then the first dynamical degree, which is birationally invariant, can be computed simply as the spectral radius of $g^*$. This also shows that rational maps which can be stabilised have a first dynamical degree that is an algebraic integer. The recent work of Bell, Diller, and Jonsson \cite{BDJ} gives an example of a rational map with transcendental first dynamical degree, hence providing another rational map which could never be stabilised.

In any case the arguments in Diller-Favre support the idea that blowing up destabilising orbits is a good approach to achieving algebraic stability.  The evidence of various examples \cite{BK06, BKTAM, BK10, Fav} shows that this approach not only succeeds in many cases but is practical to carry out explicitly.

Let us call a destabilising orbit \emph{minimal} when it does not contain any shorter ones.

\begin{prop}\label{prop:minorbit}
 Suppose that $f: X \dashto X$ is a rational map on a surface. Let\\$p, f(p), \dots, f^{n-1}(p)$ be a minimal destabilising orbit for $f$ and $\pi : X' \to X$ be the birational morphism blowing up each $p_j = f^{j-1}(p)$. Then any birational morphism $\rho : (g, Y) \to (f, X)$ stabilising $(f, X)$ factors as $\rho = \pi \circ \nu$ for some birational morphism $\nu : Y \to X'$.
 \end{prop}

\begin{defn}[\minimalalgname]\label{defn:minimal}
Given a rational surface map\\ $f_0 : X_0 \dashto X_0$ we define a (possibly finite) sequence $\pi_m : (f_{m+1}, X_{m+1}) \to (f_m, X_m)$ for $m \ge 0$ as follows
\begin{enumerate}
 \item If $f_m$ is algebraically stable, stop.
 \item If not, then pick a minimal destabilising orbit $p_1, p_2, \dots, p_n$ and blowup each of the $p_j$ to produce $\pi_m : (f_{m+1}, X_{m+1}) \to (f_m, X_m)$.
\end{enumerate}

If this sequence terminates at $f_M : X_M \dashto X_M$, write $\pi = \pi_1 \circ \cdots \circ \pi_{M-1} : X_M \to X$. Then $(f_M, X_M)$ is algebraically stable and we call $\pi : (f_M, X_M) \to (f, X)$ a \emph{minimal stabilisation} of $(f, X)$ (by blowups) when it exists.
\end{defn}

This terminology is justified by the next theorem, which says that the final result of the \minimalalgname{}, when it terminates, has a universal property.

\begin{thm}\label{thm:minimal}
 Let $f : X \dashto X$ be a rational map on a surface. If there exists a birational morphism $\rho : (g, Y) \to (f, X)$ stabilising $f$ then any instance of the \minimalalgname{} terminates in a minimal stabilisation $\pi : (\hat f, \hat X) \to (f, X)$ such that $\rho = \pi \circ \nu$ for some $\nu : Y \to \hat X$. It follows that the minimal stabilisation $(\hat f, \hat X)$ is unique for $(f, X)$.
\end{thm}

\begin{cor}\label{cor:minbirational}
 Let $f : X \dashto X$ be a birational map on a surface. Then there exists a unique minimal stabilisation $\pi : (\hat f, \hat X) \to (f, X)$.
 \end{cor}

\begin{proof}
 By \cite{DF}, there exists a birational morphism $\rho : \tilde X \to X$ which makes $\tilde f = \rho^{-1} \circ f \circ \rho$ algebraically stable on $\tilde X$. By \autoref{thm:minimal}, the minimal stabilisation $\hat f : \hat X \dashto \hat X$ via $\pi : \hat X \to X$ exists (uniquely and factors $\rho$).
\end{proof}

In certain situations, we only know how to stabilise a rational map $f$ up to a large initial iterate $f^n$, but this is perfectly sufficient for any dynamical degree computations. This weaker form of algebraic stabilisation was explored by Favre and Jonsson \cite{FJ11} for polynomial maps on $\C^2$. In their Theorem A and subsequent discussion, they show that (after possibly replacing $f$ with $f^2$) there exists a smooth surface $X$, a birational morphism $\rho : (g, X) \to (f, \P^2)$, and an integer $n \ge 1$ such that for every $m \ge 0$ \[(g^{n+m})^* = (g^n)^*(g^m)^* = (g^n)^*(g^*)^m.\] In fact, the details of their proof show that this is true for all sufficiently large $n$. Henceforth, let us say that a rational map $g : X \dashto X$ with the above property for every $n \ge N$ and $m \ge 0$ is \mbox{($N$-)}\emph{eventually algebraically stable}.
 
Likewise, we call an orbit $p, f(p)\dots, f^{m-1}(p)$ in $X$ of (closed) points \emph{$N$-eventually destabilising} iff $f^{-n}(p)$ and $f(f^{m-1}(p))$ are curves for some $n \ge N$. Note that $0$-eventual algebraic stability is the same as ordinary algebraic stability and an orbit is $0$-eventually destabilising if and only if it is a (terminal) tail of a destabilising orbit in the original sense. Our geometric criterion for algebraic stability now becomes that a rational surface map is $N$-eventually algebraically stable if and only if it has no $N$-eventually destabilising orbits.
  
  An appropriate modification of the \minimalalgname{} can be run which repeatedly blows up eventually destabilising orbits. The analogues of \autoref{prop:minorbit} and \autoref{thm:minimal} hold, giving a unique \emph{minimal eventual stabilisation} with a similar universal property; see \autoref{thm:minimalev}. The following corollary is a direct application of this to the aforementioned work by Favre and Jonsson.
  
  \begin{cor}\label{cor:minpolynomial}
 Let $f : \P^2 \dashto \P^2$ be a rational map which is a polynomial when restricted to $\C^2$. Then, after possibly replacing $f$ by $f^2$, there exists an integer $N \in \N$, and a unique minimal $N$-eventual stabilisation $\pi : (\hat f, \hat X) \to (f, X)$. More precisely, $\hat f$ is $N$-eventually algebraically stable, and if $\rho : (g, Y) \to (f, X)$ is a birational morphism where $g$ is $N$-eventually algebraically stable, then $\rho = \pi \circ \nu$ for some $\nu : Y \to \hat X$.
 \end{cor}
 
 \begin{proof}
 By (the proof of) Theorem A \cite{FJ11}, after possibly replacing $f$ with $f^2$, there exists an $N \in \N$ and a birational morphism $\rho : (g, Y) \to (f, X)$ which makes $g$ $N$-eventually algebraically stable on $Y$. By \autoref{thm:minimalev}, the \minimalalgname{} for eventual AS terminates with the unique minimal eventual stabilisation $\pi : (\hat f, \hat X) \to (f, X)$, and factors any such $\rho$.
\end{proof}

Considering only blowups over a closed point $p$ on a surface $\check X$, one can ask the same questions about algebraic stability (see \cite{GR1}). By restricting \autoref{thm:minimalev} to a neighbourhood and blowups thereof, we gain a corollary to the work on local algebraic stability by Gignac and Ruggiero \cite{GR2} who show there is an eventual algebraic stabilisation for such maps. The proof is similar to the previous corollary.

\begin{cor}\label{cor:mingerm}
  Let $(\check X, p)$ be an irreducible germ of a smooth complex surface at a point $p \in \check X$, and let $\check f : (\check X, p) \to (\check X, p)$ be a smooth morphism of germs. Let $\pi : (f, X) \to (\check f, \check X)$ be a birational morphism which blows up points only over $p$. Then, after possibly replacing $f$ by $f^2$, there exists an integer $N \in \N$, and a unique minimal $N$-eventual stabilisation $\pi : (\hat f, \hat X) \to (f, X)$ which blows up points over $\pi^{-1}(p)$.
\end{cor}

The results of \cite{DF}, \cite{Fav}, \autoref{thm:minimal}, and the above corollaries may further lead the reader to believe that if a rational map $f$ admits a stabilisation $\phi : (g, Y) \dashto (f, X)$, then in fact we can achieve algebraic stability through blowups alone, i.e. $\phi$ can be chosen to be a \emph{morphism} $\phi : (\hat f, \hat X) \to (f, X)$. This turns out to be false.

In further contrast to Gignac and Ruggiero \cite{GR2}, who show how a stabilisation through blowups can be achieved after any initial blowing up of a smooth germ, this example has no birational morphism which can stabilise it after blowing up an algebraically stable map.

\begin{thm}\label{thm:counterex} 
 Let $f : \C^2 \dashto \C^2$ be given by \[(x, y) \longmapsto (x^2, x^4y^{-3} + y^3) = \left(x^2, \frac{x^4 + y^6}{y^3}\right).\]
 Then $f$ extends to an algebraically stable rational map $f : X \dashto X$ of a Hirzebruch surface $X$. 
 If however $\sigma_0 : (f_0, X_0) \to (f, X)$ is the point blowup of $(0, 0) \in X$, then there does not exist any birational morphism $\pi : (g, Y) \to (f_0, X_0)$ which stabilises $f_0$. Furthermore this is true even if $Y$ is allowed to be singular or if we replace $f$ by an iterate $f^k$.
\end{thm}

We conclude the introduction on a somewhat different note, giving an alternative approach to the theorem of Diller-Favre for stabilising birational maps $f : X \dashto X$. In this approach one focuses on the graph of $f$, rather than on any of the destabilising orbits of $f$.

To be precise, we write $\Gamma_{\!f} \subset X \times X$ for the graph of $f$, and we call its minimal smooth desingularisation, $\sgraph_{\!f}$, the \emph{smooth graph} of $f$. We will write the maps to the first and second factor as $\alpha , \beta : \sgraph_{\!f} \to X$ respectively. Equivalently $\alpha : \sgraph_{\!f} \to X$ can be obtained as the sequence of blowups on the domain, which (minimally) resolve the indeterminacy of $f$, lifting $f$ to $\beta$; see \autoref{rmk:smoothgraph}.
  \begin{center}
 \begin{tikzcd}[column sep = 2em, row sep = 2em]
 & \sgraph_{\!f} \arrow{d}\arrow[bend right=10]{ldd}[swap, outer sep = 1.7pt]{\alpha\!}\arrow[bend left=10]{rdd}{\beta}&\\
 & \Gamma_{\!f} \arrow{ld}\arrow{rd}&\\
  X \arrow[dashed]{rr}{f} & & X
  \end{tikzcd}
\end{center}

\begin{thm}\label{thm:blowupev}
  Let $f_0 : X_0 \dashto X_0$ be a birational map on a surface.
Suppose the sequence $\alpha_m : (f_{m+1}, X_{m+1}) \to (f_m, X_m)$ is defined recursively by $X_{m+1} = \sgraph_{\!f_m}$, and $\alpha_m : \sgraph_{\!f_m} \to X_m$ is the first projection from the smooth graph. Then
\begin{enumerate}
 \item $\forall m > 0$ $X_{m+1} = \sgraph_{\!f_m} = \Gamma_{\!f_m}$, i.e.\ the graph $\Gamma_{\!f_m}$ is smooth, and
 \item $\exists M \in \N$ $\forall m \ge M$ the map $f_m$ is algebraically stable.
\end{enumerate}
\end{thm}

A key ingredient of this proof is the observation that the lift $f_1$ on $\sgraph_{\!f}$, is \emph{\nice{}}, that is whenever a curve $C$ in $\sgraph_{\!f}$ is contracted to a point by $f_1$, then $C$ contains no indeterminate points.

\begin{rmk}
 This algorithm can fail for general rational maps; see \autoref{prop:badgraphalg} for an example.
\end{rmk}

The rest of this paper is organised as follows. \autoref{sec:back} provides notation for this article and recalls useful concepts for birational maps. \autoref{sec:jud} provides the proof of \autoref{thm:minimal}. In \autoref{sec:nice} we describe some interesting properties of \nice{} birational maps and also prepare for the proof of \autoref{thm:blowupev}, which constitutes \autoref{sec:blowupev}. We end the article with examples and the computations for \autoref{thm:counterex} in \autoref{sec:counterex}.

In closing we mention that there is another recent proof of the theorem in \cite{DF} based on geometric group theory by Lonjou~and~Urech \cite{LU}. We also note that in the context of integrable systems, the failure of algebraic stability is related to the \emph{singularity confinement property}, see \cite{GRP} etc.

\section*{Acknowledgements}

I would like to thank my advisor Jeffrey Diller for his encouragement, helpful discussions, and expert advice on the exposition of this document. I thank Claudiu Raicu for his excellent question asking whether the \minimalalgname{} is truly the `minimal' method. I also thank the diligent referee, especially for the suggestion to expand the scope of this paper to the concept of eventual algebraic stability.


\section{Background}\label{sec:back}

For the rest of this article, assume all surfaces are smooth projective over an algebraically closed field, and rational maps are dominant. An account of most of the facts below can be found in \cite[\S V]{Hart}.

Let $X, Y$ be surfaces and $f : X \dashto Y$ a rational map. Let $U$ be the largest (open) set on which $f : U \to Y$ is a morphism, then we define the \emph{indeterminate set} as $I(f) = X \sm U$. Alternatively, these are the finitely many points at which $f$ cannot be continuously defined. These are often also called \emph{fundamental points}.

An irreducible curve $C \subset X$ is \emph{exceptional} iff $f(C\sm I(f))$ is a point in $Y$. We define the \emph{exceptional set}, $\E(f)$, of $f$ to be the union of all (finitely many) irreducible exceptional curves in $X$. Denote by $\comp(f)$ the number of irreducible components in $\E(f)$.

\begin{rmk}\label{rmk:rmka}
In the case of the inverse of a birational map, $I(f^{-1})$ is the proper transform of $\E(f)$, or the set of values of $f$ with infinite preimage. Using this latter definition we generalise the meaning of $I(f^{-1})$ to any rational map $f$.
\end{rmk}
 
Let $f : X \dashrightarrow Y$ be a rational map of surfaces. The \emph{graph} of $f$ is the subvariety \[\Gamma_{\!f} = \overline{\set{(x, f(x)) \in X\sm I(f) \times Y}} \subset X\times Y\] along with projections $\alpha : \Gamma_{\!f} \to X$ onto the first factor and $\beta : \Gamma_{\!f} \to Y$ onto the second factor which are proper. $\Gamma_{\!f}$ is irreducible because $X \sm I(f)$ is.

$I(f)$ is the set of points where $\alpha$ does not have a local inverse. In particular $\alpha^{-1} : X \sm I(f) \to \Gamma_f$ is an isomorphism.  For $p \nin I(f)$ we have $\beta(\alpha^{-1}(p)) = f(p)$.
 In general for any set of points $S \subseteq X$, we may define the \emph{total transform} of $S$ by $f$ as $f(S) = \beta(\alpha^{-1}(S))$. When $\emp \ne S \subseteq I(f)$ this image has dimension $1$. $\E(f) \subset X$ is the $\alpha$ projection of the set of points where $\beta$ is not finite-to-one.

$\Gamma_{\!f}$ need not be smooth and it will be more convenient to work with the minimal smooth desingularisation $\sgraph_{\!f}$ of $\Gamma_{\!f}$. Abusing notation slightly, we denote the lifted projections as $\alpha : \sgraph_{\!f} \to X$ and $\beta : \sgraph_{\!f} \to Y$. Recall the \emph{N\'eron-Severi group, $\NS(X)$}; when $Y = X$ and $f$ is birational, $\beta$ is also a birational morphism, and one can deduce that $\comp(\alpha) = \comp(\beta) = \rk\NS(\sgraph_{\!f}) - \rk\NS(X)$.

\begin{rmk}\label{rmk:smoothgraph}
 Equivalently $\sgraph_{\!f}$ and $\alpha : \sgraph_{\!f} \to X$ can be obtained as the birational morphism, blowing up $X$, which (minimally) resolves the indeterminacy of $f$. Once resolved $f : X \dashto Y$ will lift through $\alpha$ to a morphism which is precisely $\beta : \sgraph_{\!f} \to Y$. We obtain a map $\alpha \times \beta : \sgraph_{\!f} \to X\times Y$ and one can check its image is $\Gamma_f$, hence it is a smooth desingularisation which dominates the minimal one. Since the minimal desingularisation is a candidate for resolving the indeterminacy we also get a birational morphism the other direction, meaning the two definitions must agree.
\end{rmk}

\begin{rmk}[Warning]
 $\sgraph_{\!f}$ may contain curves which appear neither in the domain $X$ or the codomain $Y$, meaning both $\alpha$ and $\beta$ map the curve to a point. This issue will be rectified in \autoref{prop:curveex} and \autoref{cor:graphnice}. For an example, pick any surface with a rational curve $C$ of self-intersection $0$, and construct the birational map which blows up a point on $C$ to give an exceptional curve $D$, blows up another point of $D$ to give $E$, but then blows down $C$ followed by $D$. The second curve $D$ (i.e. its proper transform) is contracted in both the domain and range of this map, however it can be found in the smooth graph (which is the intermediate surface containing all three curves).
\end{rmk}

\begin{prop}[{See \cite[\S V 5.3 \& 5.4]{Hart}}]\label{prop:factorblowup}
 Let $f : X \to Y$ be a birational morphism of surfaces. Then $f$ can be written as a composition of $\comp(f)$ point blowups.
 
 Suppose $p \in I(f^{-1})$ and $\pi : Y' \to Y$ be the point blowup of $p$. Then $f$ factors as $\pi \circ f'$ where $f': X \to Y'$ is a birational morphism with $\comp(f') = \comp(f) -1$.
 
 Otherwise if $p \nin I(f^{-1})$ and $\rho : X' \to X$ is the point blowup of $f^{-1}(p)$. Then $f$ lifts to a birational morphism $f' = \pi^{-1} \circ f \circ \rho$ with $\comp(f') = \comp(f)$.
\end{prop}

\begin{rmk}\label{rmk:falseforrat}
 Note that \autoref{prop:factorblowup} does not hold for all rational morphisms, meaning it can fail if we replace ``birational'' with ``rational'', and consider $I(f^{-1})$ as defined in \autoref{rmk:rmka}. The second part of the proposition fails in \autoref{ex:ratfactor}.
\end{rmk}

\begin{defn}
Let $f : X \dashrightarrow Y$ and $g : Y \dashrightarrow Z$ be rational maps on surfaces.
 We say an irreducible curve $C \subset X$ is a \emph{destabilising curve} for the composition $g \circ f$ iff $f(C\sm I(f)) = y \in I(g)$. Equivalently $C \subseteq f^{-1}(y)$ and $g(y) \supseteq D$ for some irreducible curve $D \subset Z$. We call $D$ an \emph{inverse destabilising curve}, and we say that $y$ \emph{destabilises the composition} $g \circ f$. We say that the composition $g \circ f$ is \emph{locally stable at} $x \in X$ iff $x$ is not contained in any destabilising curve.
 \end{defn}
 
 By unravelling the definitions given above, we get the following proposition.
 
 \begin{prop}\label{prop:graphcomp}
 Let $f : X \dashto Y$ and $g : Y \dashto Z$ be rational maps. Then \[g(f(x)) \supseteq (g\circ f)(x)\] with equality if $x$ is not contained in a destabilising curve. Moreover $z \in g(f(x))$ if and only if $z \in (g\circ f)(x)$ or we can find $y \in I(g)\cap I(f^{-1})$ destabilising the composition $g\circ f$ where $z \in g(y)$ and $x \in f^{-1}(y)$.
\end{prop}

\begin{defn}
Let $f : X \dashto X$ be a rational map and let $p, f(p)\dots, f^{m-1}(p)$ be an orbit in $X$. 
We say this is a \emph{destabilising orbit} iff there is an (indeterminate) $q \in f^{m-1}(p)$ such that both $f(q) = D$ and $f^{-1}(p) = C$ are curves.
Let $N \in \N$, then we say the orbit is \emph{$N$-eventually destabilising} iff there is a closed point $q \in f^{m-1}(p)$ such that both $f(q) = D$ and $f^{-n}(p) = C$ are curves for some $n \ge N$.

 In either destabilising scenario, we say the \emph{length} of the orbit is $m$, we call each irreducible component of $C$ a \emph{destabilising curve} and each component of $D$ an \emph{inverse destabilising curve} of $f$.
\end{defn}

\begin{defn}
 Let $f : X \dashto X$ be a rational map. We say that $f$ is \emph{algebraically stable} iff for every $n \ge 0$ we have $(f^n)^* = (f^*)^n$. We say $f$ is \emph{$N$-eventually algebraically stable} iff for every $n \ge N$ and $m \ge 0$ we have $(f^{n+m})^* = (f^n)^*(f^m)^* = (f^n)^*(f^*)^m$.
\end{defn}

\begin{prop}[{\cite[1.14]{DF}, \cite{FS}}]
 Let $f : X \dashto X$ be a rational map. Then $f$ is algebraically stable if and only if $f$ has no destabilising orbits. Similarly, $f$ is $N$-eventually algebraically stable if and only if $f$ has no $N$-eventually destabilising orbits.
\end{prop}

 This geometric characterisation of algebraic stability can be understood through \autoref{prop:graphcomp}. One can show that $f^*g^*D - (g\circ f)^*D$ is supported on destabilising curves for any divisor $D$ and moreover that the composition $g\circ f$ is stable if and only if $f^*g^* = (g\circ f)^*$.
 
 
\section{Minimal Stabilisation}\label{sec:jud}

\begin{defn}
 We say that a destabilising orbit $p, f(p), \dots, f^{m-1}(p)$ is \emph{minimal} iff the $p_j = f^{j-1}(p)$ are all distinct closed points where for $1 \le j < m$ we have $p_j \nin I(f)$, and for $1 < j \le m$ we have $p_j \nin I(f^{-1})$.
\end{defn}

It is easy to see that the points of a minimal destabilising orbit do not contain any shorter destabilising orbits and conversely that any destabilising orbit of minimum length $m$ is \emph{minimal}, so these must always exist when $f$ is not algebraically stable. Minimality is less natural for eventually destabilising orbits. An $N$-eventually destabilising orbit always possesses a point that constitutes a singleton $N$-eventually destabilising orbit, which is vacuously minimal.

\begin{proof}[Proof of \autoref{prop:minorbit}]
 First note that by applying \autoref{prop:factorblowup} ($n$ times), if $\rho$ blows up all the $p_j$ at least once, then we get a new birational morphism $\nu$ which provides the factorisation $\rho = \pi \circ \nu$. Therefore we will proceed to show that $\rho$ does indeed blowup the $p_j$. Let $\hat p_j = \rho^{-1}(p_j)$, which a priori may be closed points or exceptional curves.

 \begin{center}
\begin{tikzcd}[column sep = 2.5em, row sep = 2.5em]
Y \arrow[swap]{d}{\rho} \arrow[dashed]{r}{g} & Y \arrow{d}{\rho}\\
 X \arrow[dashed]{r}{f} & X
\end{tikzcd}
\end{center}

Suppose not; then there is a largest $m \le n$ such that $\hat p_m$ is a closed point in $Y$, so we claim this is indeterminate for $g$. Say $m = n$ and $f(p_n) = D \subset X$, then $\hat p_n \in I(g)$ if $g(\hat p_n)$ is a curve, for which it is enough to show that $\rho(g(\hat p_n))$ is a curve. Indeed
\[\rho(g(\hat p_n)) = \rho \circ g(\hat p_n) = f \circ \rho(\hat p_n) = f(p_n) = D\] is a curve; note that $I(\rho) = \emp$ and the last step holds because $\rho$ is locally an isomorphism at $\hat p_n$. If $m < n$ then $p_{m+1} \nin I(f^{-1})$, so the composition $\rho^{-1} \circ f$ is locally stable at $p_m$, and $\rho$ is locally an isomorphism at $\hat p_m$ so $(\rho^{-1} \circ f) \circ \rho = g$ is locally stable near $\hat p_m$. Therefore \[g(\hat p_m) = \rho^{-1} (f(\rho(\hat p_m))) = \rho^{-1}(p_{m+1})\] which is also a curve by assumption.

Suppose that $k\le m$ is minimal such that $p_k, \dots, p_m$ are not blown up by $\rho$ and either $k=1$ or $k-1$ is blown up by $\rho$. Since $\rho$ is a local isomorphism over these points, one sees that $\hat p_k \mapsto \hat p_{k+1} \mapsto \cdots \mapsto \hat p_m$. To provide a contradiction to the algebraic stability of $g$, we will next show that this is a destabilising orbit for $g$. We only need to show that $g^{-1}(p_k)$ is a curve and this is very similar to the case of $p_m$ above.

$\rho$ is locally an isomorphism near $\hat p_k$. So it is enough to show that $(f \circ \rho)^{-1}(p_k)$ is a curve. Say $k=1$ and hence $f^{-1}(p_1)$ contains a curve $C$, then $(f \circ \rho)^{-1}(p_1)$ contains $\rho^{-1}(C\sm I(f))$ which is (almost all of) a curve. Otherwise $k > 1$, so $p_{k-1} \nin I(f)$ and the composition $f \circ \rho$ is algebraically stable over $\rho^{-1}(p_{k-1})$ which is a curve $\hat C$. We get that $\rho^{-1}\circ f^{-1}(p_k) = \hat C$.
\end{proof}


\begin{proof}[Proof of \autoref{thm:minimal}]
Given that $g : Y \dashto Y$ dominates $f$ via $\rho : (g, Y) \to (f, X) = (f_1, X_1)$ we may proceed inductively on $m$ with the hypothesis that $g : Y \dashto Y$ dominates $f_m : X_m \dashto X_m$ via $\nu_m : Y \to X_m$.

If $f_m$ is algebraically stable we are done, otherwise \autoref{prop:minorbit} says that because $\pi_m$ blows up a minimal destabilising orbit we have a $\nu_{m+1} : Y \to X_{m+1}$ which factors $\nu_m$ as $\nu_m = \pi_m \circ \nu_{m+1}$.
\begin{center}
\begin{tikzcd}[column sep = 2.5em, row sep = 2.5em]
 X_{m+1} \arrow[swap]{d}{\pi_m}& Y \arrow[swap, red]{l}{\nu_{m+1}} \arrow{ld}{\nu_m}\\
 X_m
\end{tikzcd}
\end{center}
Clearly there is no limit to the number of times we can do this if $f_m$ is never algebraically stable for $m \ge 1$. However overall we have shown that \[\rho = \nu_1 = \pi_1 \circ \nu_2 = \pi_1 \circ \pi_2 \circ \nu_3 = \cdots = \pi_1 \circ \cdots \circ \pi_m \circ \nu_{m+1} = \pi \circ \nu_{m+1}\] meaning that $\comp(\rho) \ge m$. Therefore $\set[m \in \N]{f_m \text{ is not AS}}$ is in fact bounded above, strictly by $m = M$ say, and whence $f_M$ is algebraically stable. Moreover, $\rho = \pi \circ \nu$ where we define $\nu = \nu_M$.

We have shown that $\pi$ factors any birational morphism stabilising $f$, and to finish we apply this to get uniqueness of $(\hat f, \hat X)$. Suppose we proceed in the \minimalalgname{} in two different ways which produce two (potentially different) models, namely $\hat f_1 : \hat X_1 \dashto \hat X_1$ via $\pi_1 : \hat X_1 \to X$ and $\hat f_2 : \hat X_2 \dashto \hat X_2$ via $\pi_2 : \hat X_2 \to X$. By the above we have that $\pi_1 = \pi_2 \circ \nu_1$ and $\pi_2 = \pi_1 \circ \nu_2$. We deduce that $\nu_1, \nu_2$ are inverse morphisms to each other, providing an isomorphism not only of surfaces but dynamical systems $\nu_1 : (\hat f_1, \hat X_1) \leftrightarrow (\hat f_2, \hat X_2)$, i.e.\ the following diagram commutes.
\begin{center}
\begin{tikzcd}[column sep = 1.5em, row sep = 2.5em]
 \hat X_1 \arrow[swap]{dr}{\pi_1} \arrow[swap, bend right=10]{rr}{\nu_1} \arrow[loop left]{}{\id} && \hat X_2 \arrow[swap, bend right=10]{ll}{\nu_2} \arrow{ld}{\pi_2} \arrow[loop right]{}{\id}\\
 &X
\end{tikzcd}
\end{center}
\end{proof}

The following analogue of \autoref{thm:minimal} for eventual algebraic stability can then be obtained by the same proof by appropriately modifying \autoref{prop:minorbit} - see \autoref{prop:minevorbit} below. In this case, the version of \minimalalgname{} for eventual algebraic stability is to blowup some minimal $N$-eventually destabilising orbit in each step of the algorithm.

\begin{thm}\label{thm:minimalev}
 Let $f : X \dashto X$ be a rational map on a surface. If there exists a birational morphism $\rho : (g, Y) \to (f, X)$ such that $g : Y \dashto Y$ is $N$-eventually algebraically stable, then any instance of the eventual version of the \minimalalgname{} terminates in a minimal $N$-eventual stabilisation $\ev\pi : (\ev f, \ev X) \to (f, X)$ such that $\rho = \ev\pi \circ \nu$ for some $\nu : Y \to \hat X$. It follows that the minimal eventually algebraically stable model $(\ev f, \ev X)$ is unique for $(f, X)$.
\end{thm}

\begin{prop}\label{prop:minevorbit}
 Let $f: X \dashto X$ be a rational map on a surface and $N \in \N$. Consider $p, f(p), \dots, f^{m-1}(p)$, a minimal $N$-eventually destabilising orbit, meaning that $f^{-n}(p)$ and $f(f^{m-1}(p))$ are curves for some $n \ge N$. Let $\pi : X' \to X$ be the birational morphism blowing up each $p_j = f^j(p)$. Suppose that $\rho : (g, Y) \to (f, X)$ is a birational morphism such that $(g, Y)$ is $N$-eventually algebraically stable. Then $\rho$ factors as $\pi \circ \nu$ for some birational morphism $\nu : Y \to X'$.
\end{prop}

\begin{proof}
As in the proof of \autoref{prop:minorbit}, we want to show that $\rho$ blows up the destabilising orbit $(p_j)$; again write $\hat p_j = \rho^{-1}(p_j)$. Suppose this is false; then there is a largest $k \le m$ such that $\hat p_k$ is a closed point in $Y$, and we first claim that this is indeterminate for $g$. The proof of this part is identical to the first half of the proof of \autoref{prop:minorbit}. In this setting, we get the following diagram.

 \begin{center}
\begin{tikzcd}[column sep = 2.5em, row sep = 2.5em]
Y \arrow[swap]{d}{\rho} \arrow[dashed]{r}{g^n} &Y \arrow[swap]{d}{\rho} \arrow[dashed]{r}{g} &Y \arrow[swap]{d}{\rho} \arrow[dotted, no head]{r} &Y \arrow[swap]{d}{\rho} \arrow[dashed]{r}{g} & Y \arrow{d}{\rho}\\
 X \arrow[dashed]{r}{f^n} & X \arrow[dashed]{r}{f} & X \arrow[dotted, no head]{r} & X \arrow[dashed]{r}{f} & X
\end{tikzcd}
\end{center}

To give a contradiction, we claim that the singleton $\set{\hat p_k}$ is an $N$-eventually destabilising orbit for $g$. Since $\hat p_k$ is indeterminate and $n+k-1 \ge n \ge N$, it is enough to show that $g^{-(n+k-1)}(\hat p_k)$ is a curve. In turn, because $\rho$ is locally an isomorphism near $\hat p_k$, it is enough to show that $(f^{n+k-1} \circ \rho)^{-1}(p_k)$ is a curve. Given such a minimal destabilising orbit, $p_j \nin I(f)$ for every $1 \le j \le k-1 < m$; therefore the $k$-fold composition $f^{n+k-1} = f \circ \cdots \circ f \circ f \circ f^n$ is stable locally along the orbit $f^{-n}(p_1) \mapsto p_1 \mapsto \cdots \mapsto p_k$. By hypothesis, $f^{-n}(p_1) = C$ is a curve. Using \autoref{prop:graphcomp} we have that $f^{n+k-1}(C\sm I(f^n)) = f(\cdots (f(f^n(C\sm I(f^n))))\cdots) = f^{k-1}(p_1) = p_k$, i.e.\ $f^{-(n+k-1)}(p_k) = C$. Also using \autoref{prop:graphcomp} we can deduce that $\rho^{-1} \circ f^{-(n+k-1)}(p_k)$ contains $\rho^{-1}(C \sm I(\rho^{-1}))$ and hence also its closure, say $\hat C$, which is a curve.
\end{proof}

\begin{rmk}\label{rmk:simultaneousstability}
 Let $f$ be an $N$-eventually algebraically stable map and $n \ge N$, then \[(f^{kn})^* = (f^n)^*(f^{(k-1)n})^* = \cdots = \left((f^n)^*\right)^k.\] Hence for any $n \ge N$ we have that $f^n$ is algebraically stable. Now suppose that $f$ is not eventually algebraically stable, but there exists a dominant model $\pi : (g, Y) \to (f, X)$ such that $g^n$ is algebraically stable for every $n \ge N$. One can check that there is a minimal model $\tilde \pi : (\tilde f, \tilde X) \to (f, X)$ with the property that $\tilde f^n$ is algebraically stable for every $n \ge N$. This is `minimal' in the sense that if $\rho : (h, Z) \to (f, X)$ is any other such model with $h^n$ AS for all $n \ge N$ then $\rho = \tilde \pi \circ \tilde \nu$ factors. The essence of the proof is repeatedly using the universal property of algebraic stability (applying \autoref{thm:minimal}); we leave the details as an exercise to the reader. In particular, if $(f, X)$ possesses a minimal $N$-eventual stabilisation $\ev\pi : (\ev f, \ev X) \to (f, X)$ then $\ev f^n$ is AS for $n \ge N$, so by considering $(g, Y) = (\ev f, \ev X)$ as above we get that $\tilde \pi : (\tilde f, \tilde X) \to (f, X)$ exists. Furthermore $\ev\pi$ factors through $\tilde \pi$ as shown in the diagram below. This begs the following question: is $\tilde \nu$ an isomorphism?
\end{rmk}

\begin{center}
\begin{tikzcd}[column sep = 1.5em, row sep = 2.5em]
 \ev X \arrow[swap]{dr}{\ev\pi} \arrow{rr}{\tilde\nu} && \tilde X \arrow{ld}{\tilde\pi}\\
 &X
\end{tikzcd}
\end{center}

\begin{qtn}
 Is a minimal $N$-eventually algebraically stable model $\ev\pi : (\ev f, \ev X) \to (f, X)$ also the minimal model of the form $\tilde \pi : (\tilde f, \tilde X) \to (f, X)$ with the property that $\tilde f^n$ is algebraically stable for every $n \ge N$?
\end{qtn}


\section{\Nice{} maps}\label{sec:nice}

\begin{defn}
 Let $f : X \dashto Y$ be a rational map. We say $f$ is \emph{\nice{}} iff \[\E(f) \cap I(f) = \emp.\]
\end{defn}

The following results are auxiliary to \autoref{thm:blowupev} but also of independent interest.

\begin{prop}\label{prop:curveex}
 Suppose $f : X \dashto Y$ is a birational map with smooth graph $\sgraph_{\!f}$ as in the diagram below. Then $f$ is \nice{} if and only if $\E(\alpha) \cap \E(\beta) = \emp$. In this case the following additional properties hold
\begin{enumerate}[ref=\emph{(\alph*)}, label=\emph{(\alph*)}]
\item $\alpha : \E(\beta) \to \E(f)$ and $\beta : \E(\alpha) \to \E(f^{-1})$ are isomorphisms;
 \item $\comp(f) = \comp(\beta)$ and $\comp(f^{-1}) = \comp(\alpha)$;
  \item $\Gamma_{\!f} \cong \sgraph_{\!f}$ is smooth.
\end{enumerate}
  \begin{center}
 \begin{tikzcd}[column sep = 2.5em, row sep = 2em]
 & \sgraph_{\!f} \arrow{d}{\pi}\arrow[bend right=20]{ldd}[swap, outer sep = 1.7pt]{\alpha\!}\arrow[bend left=20]{rdd}{\beta}&\\
 & \Gamma_{\!f} \arrow{ld}[swap]{\pi_1\hspace{-2pt}}\arrow{rd}{\hspace{-2pt}\pi_2}&\\
  X \arrow[dashed]{rr}{f} & & Y
  \end{tikzcd}
\end{center}
\end{prop}

\begin{proof}
  If $f$ is not \nice{} then there is an indeterminate point $p$ on an exceptional curve $E$. Given that $E$ is exceptional for $f$, $f(E \sm I(f)) = \beta(\alpha^{-1}(E \sm I(f))) = q$ is a point; hence there is a curve $D \subseteq \E(\beta)$ such that $\alpha(D) = E$ and $\beta(D) = q$. Since $p \in I(f)$, we have that $\alpha^{-1}(p) = C \subseteq \E(\alpha)$ is a curve and because $p \in E$ and $\alpha^{-1}(E)$ is connected, $C$ intersects $D$. This shows that $\E(\alpha) \cap \E(\beta) \ne \emp$. 
Conversely suppose $f$ is \nice{} and let $D = \beta^{-1}(q)$ be a connected component of $\E(\beta)$. Since $\sgraph_{\!f} = \sgraph_{\!f^{-1}}$ by minimality of the factorisation, we know that $\beta$ is non-isomorphic only where its image in $Y$ is indeterminate for $f^{-1}$; let $E = f^{-1}(q)$ be the exceptional curve. Observe that the composition $\alpha\circ \beta^{-1} = f^{-1}$ is stable, therefore \[\alpha(D) = \alpha(\beta^{-1}(q)) = \alpha\circ\beta^{-1}(q) = f^{-1}(q) = E.\] By \nice{}ness, $E$ contains no indeterminate points of $f$, meaning that $\alpha$ is an isomorphism in a neighbourhood of $D$. This proves both that $\E(\alpha) \cap \E(\beta) = \emp$ and the first part of (a); the second part is by symmetry. Part (b) follows from (a) and counting irreducible components. For part (c), recall that we always have a map $\pi = \alpha \times \beta : \sgraph_{\!f} \to \Gamma_{\!f} \subset X \times Y$; it is enough to show this is an isomorphism. Over $X \sm I(f)$ we have that \[\alpha : \sgraph_{\!f} \sm \E(\alpha) \to  X \sm I(f)\] is an isomorphism. On the other hand, the projection $\pi_1 : \Gamma_{\!f} \to X$ is an isomorphism over all points where the map is injective, which in this case is also $X \sm I(f)$. Hence there is an isomorphism \[\pi_1 : \Gamma_{\!f} \sm \pi_1^{-1}(I(f)) \to  X \sm I(f).\] This shows that \[\pi = \pi_1^{-1}\circ\alpha : \sgraph_{\!f} \sm \E(\alpha) \to \Gamma_{\!f} \sm \pi_1^{-1}(I(f))\] is an isomorphism when restricted to these open subsets. Similarly, considering projection to the second factor, $\pi_2 : \Gamma_{\!f} \to Y$, we deduce that \[\pi = \pi_2^{-1}\circ\beta : \sgraph_{\!f} \sm \E(\beta) \to \Gamma_{\!f} \sm \pi_2^{-1}(I(f^{-1}))\]
   is an isomorphism. Because $\E(\alpha) \cap \E(\beta) = \emp$, the sets $\sgraph_{\!f} \sm \E(\alpha)$ and $\sgraph_{\!f} \sm \E(\beta)$ form an open cover of the smooth graph, and so we have shown that $\pi$ is an isomorphism.
\end{proof}

 \begin{rmk}
The proof that the graph $\Gamma_{\!f}$ is smooth actually shows $\pi$ is isomorphic away from $\E(\alpha) \cap \E(\beta)$. Using Zariski's Main Theorem, we know that if $\pi$ is a finite birational morphism, and $\Gamma_{\!f}$ is normal, then in fact it is an isomorphism. Using Serre's criteria for normality, if a subvariety has singularities only in codimension 2 and it is locally a complete intersection inside a smooth variety, then it is normal. In this case, $\Gamma_{\!f}$ is locally a complete intersection inside the smooth variety $X \times Y$, because locally the graph has $2$ defining equations. Hence the graph is smooth if and only if $\E(\alpha) \cap \E(\beta)$ has dimension less than $1$. 
 \end{rmk}
 
 \begin{rmk}\label{rmk:untangledrat}
 When a rational map is \nice{} but not invertible the proposition is partially false. In particular the implication remains true that a map is \nice{} whenever $\E(\alpha) \cap \E(\beta) = \emp$, but the converse can fail. See \autoref{ex:untangledrat}.
 \end{rmk}

\begin{lem}\label{lem:untanglingcomp}
Let $f : X \dashto Y$ be a birational map. Suppose $f$ can be written as $h^{-1} \circ g$ where $g : X \to Z,\ h : Y \to Z$ are birational morphisms, then $f$ is \nice{}.

Moreover $\E(f) \subseteq \E(g)$ and $\E(f^{-1}) \subseteq \E(h)$ with equality if and only if the composition $h^{-1} \circ g$ has no destabilising curves. Conversely an \nice{} $f$ always has such a decomposition.
\end{lem}

\begin{proof}
We prove the first and second part by induction. We claim that if $h^{-1} \circ g$ has a destabilising curve then we can blowup $Z$ by $\pi : Z' \to Z$ to get a simpler decomposition $h'^{-1} \circ g'$ where $\comp(g') = \comp(g) - 1$ and $\comp(h') = \comp(h) - 1$. When $f = h^{-1} \circ g$ has no destabilising curves, we conclude that $f$ is \nice{}. This process terminates because $\comp(g), \comp(h)$ cannot decrease below $0$.

Suppose we have such a destabilising curve, meaning $p \in Z$ such that $g^{-1}(p)$ and $h^{-1}(p)$ are curves. Now we blowup $p \in Z$ by $\pi : Z' \to Z$ with exceptional curve $E = \pi^{-1}(p)$. By \autoref{prop:factorblowup}, $g$ factors as $\pi \circ g'$ with $\comp(g') = \comp(g) - 1$ and $h$ factors as $\pi \circ h'$ with $\comp(h') = \comp(h) - 1$.

    \begin{center}
 \begin{tikzcd}[column sep = 1.5em, row sep = 2.5em]
  C \arrow[draw=none]{r}[sloped,auto=false]{\subseteq} & X \arrow{rrd}[swap]{g'} \arrow[bend right]{rrdd}[swap]{g} \arrow[dashed]{rrrr}{f}&&&& Y \arrow{lld}{h'} \arrow[bend left]{lldd}{h} &  D \arrow[draw=none]{l}[sloped,auto=false]{\supseteq} \\
&& E\arrow[draw=none]{r}[sloped,auto=false]{\subseteq} &Z' \arrow{d}{\pi} & ~\\
  && p \arrow[draw=none]{r}[sloped,auto=false]{\in}& Z
  \end{tikzcd}
  \end{center}
 
The proper transform of $E$, $C = \overline{g'^{-1}(E \sm I(g'^{-1}))} \subset X$, is an irreducible curve; or more simply we have $g'(C) = E$. Similarly there is an irreducible curve $D$ such that $h'(D) = E$. Whence $D$ is the proper transform of $C$ by $f = h'^{-1} \circ g'$. This completes the claim.
 
 It remains to show that if the decomposition $f = h^{-1} \circ g$ has no destabilising curves then $f$ is \nice{}, plus $\E(g) = \E(f)$ and $\E(h) = \E(f^{-1})$. Let $C \subseteq \E(g)$ be a curve, then $g(C) = p \in Z$ is a closed point. If $p \nin I(h^{-1})$ then clearly $f = h^{-1} \circ g$ is continuous on $C$ with a closed point as an image, therefore $I(f) \cap C = \emp$ and $C \subseteq \E(f)$. Otherwise if $p$ is indeterminate, then we know there is a curve $D \subseteq h^{-1}(p) \subseteq \E(h)$; whence $C$ is a destabilising curve which pairs with the inverse destabilising curve $D$, contradicting our assumption. Thus $\E(g) \subseteq \E(f)$. Conversely, if $f(C \sm I(f))$ is a closed point then certainly $g(C \sm I(f))$ is a closed point since $\E(h^{-1}) = \emp$; so $C \subseteq \E(g)$. We have shown that $\E(g) = \E(f)$; a similar argument shows that $\E(h) = \E(f^{-1})$.
 
 The converse can be seen from two perspectives. Firstly, $\alpha : \E(\beta) \to \E(f)$ is an isomorphism, hence the blowdown of $\E(\beta)$ performed by $\beta$ can be transferred (locally) to a blowdown $g$ of $\E(f)$ to an equally smooth surface $Z$. Untangledness also shows that the induced map $h : Y \to Z$ is a birational morphism. The more rigourous way to construct $Z$ is by two `charts': $U_1 = X \sm \E(f)$ and $U_2 = Y \sm \E(f^{-1})$. The transition map is $f : X \sm (\E(f) \cup I(f)) \to Y \sm (\E(f^{-1}) \cup I(f^{-1}))$ which is clearly an isomorphism. The map $g : X \to Z$ is given by gluing $\id : X \sm \E(f) \to U_1$ and $\id \circ f : X \sm I(f) \to Y \sm \E(f^{-1}) \to U_2$. On the overlap, these germs differ by exactly the transition function $f$.
\end{proof}

\begin{cor}\label{cor:graphnice}
  Let $f : X \dashto X$ be a birational map and let $\hat f$ be the lift of $f$ defined by the following commutative diagram, where $\sgraph_{\!f}$ is the smooth graph of $f$. Then $\hat f$ is \nice{}.
    \begin{center}
 \begin{tikzcd}[column sep = 1.5em, row sep = 2.5em]
  & \sgraph_{\!f} \arrow{ld}[swap]{\alpha} \arrow[dashed]{rr}{\hat f} \arrow{rd}{\beta}&& \sgraph_{\!f} \arrow{ld}{\alpha} \\
  X \arrow[dashed]{rr}{f} && X
  \end{tikzcd}
  \end{center}
\end{cor}

\begin{proof}
 $\hat f$ is \nice{} due to \autoref{lem:untanglingcomp} because $\hat f = \alpha^{-1} \circ \beta$.
 \end{proof}

\begin{cor}\label{cor:equivexceptionals}
  Let $f : X \dashto X$ be an \nice{} birational map and let $\hat f$ be the lift of $f$ defined as above. Then $\alpha : \E(\hat f) \hookrightarrow \E(f)$ and $\beta : \E(\hat f ^{-1}) \hookrightarrow \E(f ^{-1})$ are injections, which are surjective if and only if $f$ has no length $1$ destabilising orbits.
\end{cor}

\begin{proof}
 By \autoref{lem:untanglingcomp}, we know that $\E(\hat f) \subseteq \E(\beta)$ and $\E(\hat f^{-1}) \subseteq \E(\alpha)$ with equality if and only if $\alpha^{-1} \circ \beta$ has no destabilising curves; on the other hand \autoref{prop:curveex} says that $\alpha : \E(\beta) \to \E(f)$ and $\beta : \E(\alpha) \to \E(f^{-1})$ are isomorphisms. Therefore to conclude we only need to show that the composition $\alpha^{-1} \circ \beta$ has a destabilising curve $\hat C$ and inverse destabilising curve $\hat D$ if and only if the composition $f\circ f$ has a destabilising curve $C$ and inverse destabilising curve $D$. Indeed, by the isomorphisms in \autoref{prop:curveex}, $C = \alpha(\hat C)$ is a curve and $f(C) = p$ if and only if $\hat C$ is a curve and $\beta(\hat C) = p$. Similarly $D= \beta(\hat D)$ is a curve and $f^{-1}(D) = p$ if and only if $\hat D$ is a curve and $\alpha(\hat D) = p$.
\end{proof}

 
\section{Stabilisation Through Graphs}\label{sec:blowupev}
 
\begin{defn}
 Let $\mathcal D(X, f)$ be the set of all triples $(C, D, n)$ such that $C$ is a destabilising curve for $f : X \dashto X$ with an orbit of length $n$ and inverse destabilising curve $D$.
 \end{defn}

The following proposition is the heart of \autoref{thm:blowupev}.
 
\begin{prop}\label{prop:blowupev}
 Let $f : X \dashto X$ be an \nice{} birational map and $\hat f : \sgraph_{\!f} \dashto \sgraph_{\!f}$ be the lift described above in \autoref{cor:graphnice}.
Then there exists a well defined injection
\begin{align*}
\iota : \mathcal D(\sgraph_{\!f}, \hat f) &\longrightarrow \mathcal D(X, f) \\
(\hat C, \hat D, n) &\longmapsto (\alpha(\hat C), \beta(\hat D), n+1).
\end{align*}
If $\iota$ is surjective (a bijection) then $\comp(\hat f) = \comp(f)$.
If $\iota$ isn't surjective then $\comp(\hat f) < \comp(f)$.
\end{prop}

\begin{proof}
To justify that $\iota$ is well defined, we claim that every destabilising orbit upstairs descends to one downstairs; to be precise, if $(\hat C, \hat D, n) \in \mathcal D(\sgraph, \hat f)$ then $(C, D, n+1) \in \mathcal D(X, f)$ where $C = \alpha(\hat C)$ and $D = \beta(\hat D)$.

Assume $\hat C$ is a destabilising curve, $\hat f(\hat C) = \hat p$, $\hat f^{n-1}(\hat p) \ni \hat q \in I(\hat f)$ and $\hat f(\hat q) \supseteq \hat D$. Let $D = \beta(\hat D),\ C = \alpha(\hat C)$; as shown in \autoref{cor:equivexceptionals}, $C \subseteq \E(f)$ with $f(C) = p = \alpha(\hat p)$, and $D \subseteq \E(f^{-1})$ with $q = \beta(\hat q) = f^{-1}(D)$. To complete the claim we will show that $f^n(p) \ni q$.
 Indeed, consider the composition $f^n = \beta \circ \hat f^{n-1} \circ \alpha^{-1}$, which is stable by \autoref{prop:graphcomp} because $\E(\alpha^{-1}) = \emp$ and $I(\beta) = \emp$; note also that $\alpha^{-1}(p) \ni \hat p$ and $\beta(\hat q) = q$.
 \[f^n(p) = \beta ( \hat f^{n-1} ( \alpha^{-1}(p))) \supseteq \beta ( \hat f^{n-1} (\{\hat p\})) \supseteq \beta(\{\hat q\}) = \{q\}\]
 Injectivity follows from the injectivity given in \autoref{cor:equivexceptionals} and the simple fact that $\iota(\hat C, \hat D, m) = (C, D, n)$ implies $m=n - 1$.
 
For the surjectivity, we claim that $\iota$ is surjective if and only if we cannot find a length $1$ destabilising orbit for $f$, that is $(C, D, 1) \in \mathcal D(X, f)$. Then \autoref{cor:equivexceptionals} finishes the proof since we know $\comp(\hat f) \le \comp(f)$ with equality if and only if $f$ has no length $1$ destabilising orbits.

Clearly, if $(C, D, 1) \in \mathcal D(X, f)$ then $(C, D, 1)$ cannot have a preimage under $\iota$ since no destabilising orbit has length $1-1=0$. Conversely we show in the remainder of this proof that when $\mathcal D(X, f)$ has no such triples we can find a preimage for $(C, D, n) \in \mathcal D(X, f)$. Write $\hat C = \alpha^{-1}(C)$, $\hat D = \beta^{-1}(D)$, $f(C) = p$, and $q = f^{-1}(D)$. Because $f$ has no length 1 destabilising orbits we know that $p \nin I(f) = I(\alpha^{-1})$ and $q \nin I(f^{-1}) = I(\beta^{-1})$; furthermore by \autoref{cor:equivexceptionals}, $\hat C \subseteq \E(\hat f)$, $\hat D \subseteq \E(\hat f^{-1})$ are both irreducible curves. Hence on $\sgraph_{\!f}$ we have two closed points $\hat p = \hat f(\hat C) = \alpha^{-1}(p)$ and $\hat q = \hat f^{-1}(\hat D) = \beta^{-1}(q)$. To complete the destabilising orbit between $\hat C$ and $\hat D$, we wish to show $\hat q \in \hat f^{n-2}(\hat p)$ given that $q \in f^{n-1}(p)$. Considering $p \nin I(\alpha^{-1})$, $q \nin I(\beta^{-1})$ we may apply \autoref{prop:graphcomp} twice to obtain
 \[\hat f^{n-2}(\hat p) = \beta^{-1} \circ f^{n-1} ( \alpha(\hat p)) = \beta^{-1} \circ f^{n-1} (p) \supseteq \beta^{-1}(\{q\}) = \{\hat q\}.\]
 Therefore $(\hat C, \hat D, n-1) \in \mathcal D(\sgraph_{\!f}, \hat f)$ and $\iota((\hat C, \hat D, n-1)) = (C, D, n)$.
  \end{proof}

\begin{proof}[Proof of \autoref{thm:blowupev}]
 First, note that by \autoref{prop:curveex} $X_m$ is smooth for all $m \ge 1$, then by \autoref{cor:graphnice}, $f_m$ is \nice{} for $m \ge 1$; assume that $(X_1, f_1)$ is not algebraically stable.
 
 If $f_m$ is not algebraically stable then we may choose $(C, D, n) \in \mathcal D(X_m, f_m)$. By \autoref{prop:blowupev}, $\comp(f_{m+1}) \le \comp(f_m)$ and either we have strict inequality or all destabilising orbits lift to strictly shorter destabilising orbits. Since lengths of orbits must be positive, eventually we find an $m' \le m+n$ such that $\comp(f_{m'}) < \comp(f_m)$.
 
 The sequence $\comp(f_m) \ge 0$ must stabilise as $m$ increases with $\comp(f_m) = \comp(f_M)$ for all $m \ge M$. Then $\mathcal D(X_m, f_m) = \emp$ for all such $m$, otherwise we could decrease $\comp(f_m)$ further as above.
\end{proof}

\autoref{prop:blowupev} and the theorem admit a simpler proof, without counting $\comp(f_m)$, when $\mathcal D(X, f)$ is a finite set or the lengths of destabilising orbits are bounded. One can show that a length $n$ destabilising orbit on $\sgraph_{\!f}$ implies the existence of a length $n+1$ orbit on $X$. The key is to note that for any birational $f$, it would be absurd to have any connected component of $\E(\alpha)$ contracted by $\beta$. Therefore, in any connected component of the exceptional locus $\E(\hat f) \subseteq \E(\beta)$ we can find an irreducible curve $\hat C \subset \E(\hat f)$ which does \emph{not} get contracted by $\alpha$. The rest is similar to the first two paragraphs of the proof of \autoref{prop:blowupev}. This gives the following quantitative bound.

\begin{cor}
Let $f : X \dashto X$ be a birational map and the sequence $(X_m, f_m)$ be as given in \autoref{thm:blowupev}. Suppose that $N$ is an upper bound on lengths of a destabilising orbit for $f$. Then $\forall m \ge N$, $f_m$ is algebraically stable.
\end{cor}

\section{Examples}\label{sec:counterex}

\subsection{Introduction}

The main goal of this section is to prove \autoref{thm:counterex}; later we apply the same techniques to other examples. 
This demonstrates that there are rational maps which can be stabilised by a birational conjugacy but not by blowups of the surface alone. 
Recall that $f : \C^2 \dashto \C^2$ is given by \[(x, y) \longmapsto (x^2, x^4y^{-3} + y^3) = \left(x^2, \frac{x^4 + y^6}{y^3}\right).\]

Let us initially compactify $(f, \C^2) \dashto (\tilde f, \P^1 \times \P^1)$ in the obvious way, by adding a point at infinity to each factor. Unfortunately, $\set{x=\infty}$ is a destabilising curve for $\tilde f$. One may check that $\tilde f(\set{x=\infty})$ is the indeterminate point $(\infty, \infty)$. However if we modify this fibre to create $X$, the first Hirzebruch surface, then $f : X \dashto X$ is algebraically stable.

\begin{prop}\label{prop:initialstability}
 Let $\psi : (\tilde f, \P^1 \times \P^1) \dashto (f, X)$ be the birational transformation obtained by blowing up $(\infty, \infty)$ and then blowing down the proper transform of $\set{x = \infty}$. Then the exceptional set, $\E(f)$, is $E_\infty = \psi(\infty, \infty)$, and $f(E_\infty) \nin I(f)$ is a fixed point. In particular, $f$ is algebraically stable.
\end{prop}

This proposition can be verified by local coordinate computations, or by following the method laid out in \S\ref{sec:mapexc}. 

Since $\psi$ acts trivially on $\C \times \P^1$ we will think of this as a subset of $X$; we define $E_0 = E_{\frac 01}$ to be $\set{x=0} \subset \C \times \P^1$. Next we define the blowup $\sigma_0 : (f_0, X_0) \to (f, X)$ centred on $(0, 0) \in \C^2 \subset X$ and let $E_1 = E_{\frac 11} = \E(\sigma_0)$ be the exceptional curve. A consequence of our proof will be that $f_0$ is not algebraically stable, but we will also see this directly from \autoref{prop:mapexc}, later.

The method we provide in the next two subsections is elementary, however analysing these examples on Berkovich space is much faster and more informative. The reader who is familiar with Berkovich theory may skip to \S\ref{sec:berkalt}.

\subsection{Satellite Blowups}

Before analysing $f$ further, we introduce a convenient bookkeeping system for blowups over the origin. See \cite[\S 6.1]{FJ04}, \cite[\S 15.1]{Jon} and \cite[appx]{DL} for precedents.

\begin{defn}
 A birational morphism $\phi : Y \to X_0$ is \emph{satellite (relative to $E_0$ and $E_1$)} iff $\phi = \sigma_1 \circ \sigma_2 \circ \cdots \circ \sigma_n$ ($n \ge 1$) and for every $1 \le j \le n$ we have that $\sigma_j$ is the blowup of the point of intersection between two curves in the list $E_0, E_1, \E(\sigma_1), \E(\sigma_2), \dots , \E(\sigma_{j-1})$.
\end{defn}

Note that in this definition and throughout this section we will adopt the convention that if $\phi : Y \to Z$ is a birational morphism between surfaces and $C \subset Z$ is a curve then $C$ will also denote the proper transform $\overline{\phi^{-1}(C)\sm \E(\phi)}$ of $C$ in $Y$.

 Now given any birational morphism $\phi : Y \to X_0$ which is \emph{satellite}, as above, we proceed to index the exceptional curves of $\phi$, by rational numbers $\frac ab \in (0, 1)$ in lowest terms as follows. For each $j \ge 1$, if $\sigma_j$ blows up the intersection $E_{\frac ab} \cap E_{\frac cd}$ of two previously indexed curves from among $E_{\frac01}, E_{\frac11}, \E(\sigma_1), \E(\sigma_2), \dots , \E(\sigma_{j-1})$, then we declare $\E(\sigma_j) = E_{\frac{a+c}{b+d}}$. Note that the Farey sum $\frac{a+c}{b+d} \in (\frac ab, \frac cd)$ is a rational number in lowest terms.

Let $0 = r_0 < r_1 < \cdots < r_n = 1$ be the full list of rational indices for the curves $E_{r_j}$ as above. The \emph{dual graph} for $\phi$ is defined as the graph with vertices $\set[E_{r_j}]{0\le j \le n}$ and edges $\set[E_{r_j}E_{r_k}]{E_{r_j} \pitchfork E_{r_k}}$. 
\begin{center}
 \begin{tikzcd}
 E_{r_0} \arrow[-]{r} & E_{r_1} \arrow[-]{r} & \cdots \arrow[-]{r} & E_{r_{n-1}} \arrow[-]{r} & E_{r_n}
 \end{tikzcd}
\end{center}
In particular blowing up $E_{\frac ab} \cap E_{\frac cd}$ corresponds to inserting a vertex as follows.
\begin{center}
 \begin{tikzcd}
 \cdots \arrow[-]{r} & E_{\frac ab} \arrow[-]{rr} &\arrow[rightsquigarrow]{d}& E_{\frac cd} \arrow[-]{r} & \cdots\\
 \cdots \arrow[-]{r} & E_{\frac ab} \arrow[-]{r} & E_{\frac{a+c}{b+d}} \arrow[-]{r} & E_{\frac cd} \arrow[-]{r} & \cdots
 \end{tikzcd}
\end{center}
We caution however that the ordering of the dual graph does not match the ordering of the $\sigma_j$, i.e.\ $E_{r_j} \ne \E(\sigma_j)$ in general.

The curves $E_{\frac ab}$ can be seen as `degenerations' of embeddings of the complex torus $\C^*\times \C^*$ into $Y$. For any $\frac ab$ we define the map $\gamma_{\frac ab} : \P^1 \times \P^1 \dashto Y$ given on $\C^*\times \C^*$ by
\begin{align*}
 \gamma_{\frac ab} : \C^*\times \C^* &\longrightarrow \C^*\times \C^* \subset Y\\
 (s, t) &\longmapsto (t^b, st^a).
\end{align*}
In fact this map represents the toric blowup over the origin with weight $(a, b)$. When $Y$ is such a toric blowup, $\gamma_{\frac ab}$ gives coordinates $(s, t)$ for a particular chart at the exceptional curve of weight $(a, b)$. Indeed, \autoref{prop:degen} states that the locus $\set{t=0}$ corresponds to the exceptional curve $E_{\frac ab}$.

More generally, we say that a rational map $\gamma : \P^1 \times \P^1 \dashto Y$ is asymptotic to $\gamma_{\frac ab}$, or $\gamma \sim \gamma_{\frac ab}$, if and only if on $\C^* \times \C^*$ we have
 \[(s, t) \longmapsto (t^d + \lo(t^d), R(s)t^c + \lo(t^c))\]
where $\frac cd = \frac ab$ and $R$ is a non-constant rational function on $\P^1$.

\begin{prop}\label{prop:degen}
Suppose that $\gamma \sim \gamma_{\frac ab} : \P^1\times \P^1 \dashto Y$. Let $Z$ be the proper transform of $\P^1 \times \set 0$ under $\gamma$. Then either
\begin{enumerate}
 \item \label{prop:degen:a} $\frac ab = r_j$ for some $j$, and $Z$ is the curve $E_{\frac ab}$; or
 \item $r_j < \frac ab < r_{j+1}$ and $Z$ is the closed point $E_{r_j}\cap E_{r_{j+1}}$.
\end{enumerate}
\end{prop}

Part \ref{prop:degen:a} can be proven inductively by resolving the singularities of $\gamma_\frac ab$, starting with $\frac ab \in \Z$ and using Farey addition. If $Z$ is not a curve in a particular model $Y$, then it must be a closed point. To show that $r_j < \frac ab < r_{j+1}$ is analogous to proving on $Y = \P^1 \times \P^1$ that $Z = \set{(0, 0)}$ for any $\frac ab > 0 = r_0$.

\subsection{Mapping Exceptional Curves}\label{sec:mapexc}
 
Now we compute the images of the curves $E_{r_j}$ under $g$, using \autoref{prop:degen}. The rough idea is that the point or curve represented by $\gamma_{\frac ab}$ is mapped to the point or curve represented by $f \circ \gamma_{\frac ab} \sim \gamma_{\frac cd}$. 
In \autoref{prop:mapexc} we will state for example how $E_{r_j}$ maps to another curve $E_{r_k}$ or is contracted to a point $E_{r_k} \cap E_{r_{k+1}}$, depending on whether $\frac cd = r_k$ or $\frac cd \in (r_k, r_{k+1})$ respectively. The rest of the proposition describes the image of a (possibly indeterminate) point $E_{r_j} \cap E_{r_{j+1}}$. We proceed to compute $\frac cd$ in terms of $\frac ab$.
 \[f\circ \gamma_{\frac ab}(s, t) = f(t^b, st^a) = \left(t^{2b}, s^{-3}t^{4b-3a} + s^3t^{3a}\right)\]
In the case where $4b-3a > 3a$, looking at lowest order terms, we get that $f\circ \gamma_{\frac ab} \sim \gamma_{\frac {3a}{2b}}$. In the case where $4b-3a < 3a$ we get that $f\circ \gamma_{\frac ab} \sim \gamma_{\frac {4b-3a}{2b}}$. Finally in the special case that $4b-3a = 3a$ we get $f\circ \gamma_{\frac ab}(s, t) = \left(t^{2b}, (s^{-3}+ s^3)t^{3a}\right) \sim \gamma_{\frac {3a}{2b}}(s, t)$. 
In short $f \circ \gamma_q \sim \gamma_{T_f(q)}$, where
\[T_f : q\longmapsto 
\begin{cases}
 \frac {3}{2}q & q \le \frac 23\\
 2 - \frac {3}{2}q & q > \frac 23.
\end{cases}
\]

\begin{ex}
 We can now see that $f_0 : X_0 \dashto X_0$ is not algebraically stable. $E_0$ is fixed by $f_0$ since $T_f(0) = 0$. However $T_f(1) = \half \in (0, 1) = (r_0, r_1)$, therefore $f_0(E_1) = E_0 \cap E_1 = P$. This point is indeterminate with $f_0(P) = E_1$ because $T_f((0, 1)) = (0, 1] \ni 1$.
\end{ex}

\begin{prop}\label{prop:mapexc}
 Let $\phi : (g, Y) \to (f_0, X_0)$ be satellite as above, and let the irreducible curves of the fibre $\set{x=0}$ on $Y$ be indexed by rational (Farey) parameters \[0 = r_0 < r_1 < \cdots < r_n = 1.\] Let $T_f : \Q \cap [0, 1] \to \Q \cap [0, 1]$ be such that $g \circ \gamma_q = \gamma \sim \gamma_{T_f(q)}$ for every $q \in \Q \cap [0, 1]$.
 
 Then the dynamics over $\set{x=0}$ is determined by the following:
\begin{enumerate}
 \item if $q = r_j$ and $T_f(q) = r_k$ for some $0 \le j, k \le n$, then $g$ maps $E_q$ to $E_{T_f(q)}$;
 \item if $q = r_j$ and $r_k < T_f(q) < r_{k+1}$ for some $0 \le j, k \le n$, then $g : E_q \mapsto E_{r_k} \cap E_{r_{k+1}}$;
 \item if $T_f((r_j, r_{j+1})) \subset (r_k, r_{k+1})$, then $g\left(E_{r_j} \cap E_{r_{j+1}}\right) = E_{r_k} \cap E_{r_{k+1}}$; otherwise
 \item if $[r_k, r_l] \subseteq T_f((r_j, r_{j+1}))$ with $k$ minimal and $l$ maximal, then we have \[g\!\left(E_{r_j} \cap E_{r_{j+1}}\right) = E_{r_k} \cup \cdots \cup E_{r_l}.\]
\end{enumerate}
\end{prop}

\subsection{Dynamics of Exceptional Curves}\label{sec:dynexc}

 Suppose for contradiction we have a birational morphism $\pi : (g, Y) \to (f_0, X_0)$ which stabilises $f$. Then by \autoref{thm:minimal} we may assume that $\pi : (g, Y) = (\hat f, \hat X) \to (f_0, X_0)$ is the minimal stabilisation. Consider precisely how $\pi$ blows up $X_0$ with the \minimalalgname{}.

\begin{lem}\label{lem:satellitestab}
 The \minimalalgname{} on $(f_0, X_0)$ only creates curves which are satellite relative to $E_0$ and $E_1$.
\end{lem}

Hence we can assume that $\E(\pi) = E_{r_1} \cup \cdots \cup E_{r_{n-1}}$, $r_j \in (0, 1) \cap \Q$. Note that the interval $[0, 1]$ is forward invariant for $T_f$, as is $[\half, 1]$. The rest of the proof will hinge on $T_f$ being topologically mixing on this interval. In the Berkovich theory, this lemma also corresponds to the similar fact that $[\zeta(0, 1), \zeta(0, |x|)]$ is forward invariant; see \autoref{sec:berkalt} for further details.

\begin{proof}
The only possible destabilising curve for $f_0$ is $E_1$. Any destabilising orbit beginning at $E_1$ must remain in the fibre $\set{x=0}$ which is fixed by $f_0$, and the same applies for every exceptional curve of $\pi$ created by the algorithm.

 At the first step of the algorithm we only have $E_0, E_1$. Proceeding inductively, suppose we have an intermediate surface $\pi' : (g, Y) \to (f_0, X_0)$ generated by the algorithm which is satellite relative to $E_0$ and $E_1$, with $\E(\pi') = E_{r_1'} \cup \cdots \cup E_{r_m'}$. 
 On one hand, a destabilising curve must be one of the $E_{r_j'}$, but on the other hand, \autoref{prop:mapexc} says that a minimal destabilising orbit consists of finitely many points of the form $E_{r_j'}\cap E_{r_{j+1}'}$. Blowing up all of these points leads to a further map which is satellite relative to $E_0$ and $E_1$.
\end{proof}

\begin{proof}[Proof of \autoref{thm:counterex}]
Consider the dynamics of the expanding piecewise affine map $T_f$.
\[1 \mapsto \frac 12 \mapsto \frac 34 \mapsto \frac 78 \mapsto \frac{11}{16} \mapsto \cdots\]
Suppose that $\frac ab = \frac a{2^n}$ with $a$ odd, then
\[\frac a{2^n} \longmapsto
\begin{cases}
 \frac {3a}{2^{n+1}} & \mathmakebox[4mm]{\frac a{2^n}} < \frac 23\\
 2 - \frac {3a}{2^{n+1}} = \frac{2^{n+2} - 3a}{2^{n+1}} & \mathmakebox[4mm]{\frac ab} > \frac 23
\end{cases}
\]
where both $3a$ and $2^{n+2} -3a$ are also odd. 
In particular, the we see that $(T_f^m(1))$, the orbit of $1$, is infinite. Hence there exists a smallest $m = m_1$ such that $T_f^{m-1}(1) = q \in \set{r_0, \dots, r_n}$, but $T_f(q) \in (r_j, r_{j+1})$ for some $j$. Thus by \autoref{prop:mapexc}(i, ii), $\hat f(E_q)$ is the closed point $P = E_{r_j} \cap E_{r_{j+1}}$. We claim that for some (smallest) $m = m_2$, the interval $T_f^m((r_j, r_{j+1}))$ contains one of the $r_l$. Then by \autoref{prop:mapexc}(iii, iv) we have that for $P, \hat f(P), \dots, \hat f^{m_2-1}(P) = Q$ are closed points with $\hat f(Q) \supset E_l$. Therefore $\hat f : \hat X \dashto \hat X$ still has the (minimal) destabilising orbit $P, \hat f(P), \dots, \hat f^{m_2-1}(P) = Q$, contradiction.

 The map $T_f(q)$ is expanding lengths by a factor of $\frac 32$ on all intervals which do not include $q = \frac 23$. Since $T_f([0, 1]) = [0, 1]$, this cannot occur indefinitely, therefore one of the intervals $T_f^m((r_j, r_{j+1}))$ contains $\frac 23$, and since $T_f(\frac 23) = 1 = r_n$, we have proved the claim. We can actually go much further, having shown that eventually $1 \in T_f^m((r_j, r_{j+1}))$, we can deduce that there is an $M \in \N$ such that for all $m \ge M$, $T_f^m((r_j, r_{j+1})) = [\half, 1]$.
 
To see that even $f^k$ cannot be stabilised for any $k \in \N$, we instead consider the corresponding map on indices $T_f^k = T_{f^k}$, but the same argument works. Again we have that $(T_f^{km}(1))$ is infinite, and by the extended claim above, eventually $T_f^{km}((r_j, r_{j+1})) \ni r_l$ for some $l$.

Finally we justify the theorem in the case of arbitrary (singular) stabilisations. Consider a birational morphism $\pi : (g, Y) \to (f_0, X_0)$ with $Y$ singular. First, we can factor $\pi$ as $\pi' \circ \phi^{-1}$ where $\pi' : (g', Y') \to (f_0, X_0)$ and $\phi : (g', Y') \to (g, Y)$ are birational morphisms, $Y'$ is smooth, and $\phi$ contracts only some curves in $\E(\pi')$ (whose self intersection is negative). Second, if $\pi'$ has a non-satellite component, then one can check that $\pi'$ can be factored as $\pi_s \circ \pi_f$ where $\pi_s : (h, Z) \to (f_0, X_0)$ and $\pi_f : (g', Y') \to (h, Z)$ are birational morphisms such that $\pi_s$ is satellite between $E_0$ and $E_1$, and $\pi_f$ is a birational morphism that is isomorphic over a neighbourhood of every $E_q\cap E_r$. By the first part of the proof, $h : Z \dashto Z$ has a destabilising orbit $P, \dots, Q$ along points of the form $E_{r_j}\cap E_{r_{j+1}}$. Since $\pi_f : Y' \to Z$ is isomorphic over these points, the destabilising orbit persists for $g' : Y' \dashto Y'$. Moreover, the first part of the proof shows that we can assume the destabilising curve and inverse destabilising curve are $E_1$. More precisely, for some (non-minimal) $m_1, m_2$, we have $E_1 \subseteq (g')^{-m_1}(P)$, and $E_1 \subseteq (g')^{m_2}(Q)$. Note that $\phi$ may contract some of the $E_{r_j}$ from our list of divisors in $Y'$, but it preserves at least $E_1$ since it exists in $X_0$. Therefore after projecting by $\phi$, the destabilising orbit persists for $g : Y \dashto Y$.
\end{proof}

\subsection{The Berkovich Alternative}\label{sec:berkalt}
Here we provide some details about another approach to the bookkeeping using the Berkovich projective line, $\P^1_\text{an}(K)$. The following notation can be found in \cite{Bene}. We also refer the reader to the theory of `skew products' on the Berkovich projective line, as developed in the author's PhD thesis \cite{thesis}. We work over the field, $K$, of Puiseux series in the variable $x$ (the same $x$ as above) with $\C$ coefficients and the $x$-adic norm. When $q \in \Q\cap[0, 1]$, there is a correspondence between Type II norms of the form $\zeta(0, |x|^q)$, disks of Puiseux series $\overline D(0, |x|^q)$, and the exceptional curves $E_q$ which can be obtained from satellite blowups between $E_0$ and $E_1$. The family of curves $t \mapsto \gamma_\frac ab(s, t)$ parametrised by $s$ is precisely giving us a dense family of Puiseux series $y = sx^{\frac ab}$ in $\overline D(0, |x|^{\frac ab})$. A Type II norm $\zeta(0, |x|^{\frac ab})$ can be obtained as a maximum over a dense family of Type I seminorms, such as $sx^{\frac ab} \in \P^1(K) \subset \P^1_\text{an}(K)$, in the corresponding disk. This Type II norm measures the order of vanishing of functions at the particular prime divisor over $\set{x=0}$ where the curves $y = sx^{\frac ab}$ land, namely $E_{\frac ab}$. Hence the curves $E_{r_j}$ correspond to finitely many Type II points $\zeta(0, |x|^{r_j})$ in $[\zeta(0, 1), \zeta(0, |x|)]$, and an intersection point $E_{r_j}\cap E_{r_{j+1}}$ corresponds to the Berkovich annulus $U_j = \set{|x|^{r_{j+1}} < |\zeta| < |x|^{r_j}}$ bounded by $\zeta(0, |x|^{r_j})$ and $\zeta(0, |x|^{r_{j+1}})$.

The second component of our rational map induces a rational map $y \mapsto x^4y^{-3} + y^3$ on $\P^1_\text{an}$. On $(0, \infty) \subset \P^1_\text{an}$ this maps $\zeta(0, r)$ to $\zeta(0, R)$, where $R$ is the radius given by the magnitude of the Laurent series $x^4y^{-3} + y^3$ at $|y| = r$. The Weierstrass degree is $-3$ when $R = |x|^4r^{-3} = |x^4y^{-3}| > |y^3| = r^3$ and $3$ when $|x|^4r^{-3} = |x^4y^{-3}| < |y^3| = r^3 = R$. This means that
 \[\zeta(0, r) \longmapsto 
\begin{cases}
 \zeta(0, r^3) & r > |x|^{\frac 23}\\
 \zeta(0, |x|^4r^{-3}) & r < |x|^{\frac 23}.
\end{cases}
\]
The effect of the first component of $f$, $x \mapsto x^2$ is to replace each diameter with its square root. The map $T_f$ constructed in \S\ref{sec:mapexc} describes the dynamics on $(0, \infty) \subset \P^1_\text{an}$ with each $T_f(q) = r$ corresponding to $f_*(\zeta(0, |x|^q)) = \zeta(0, |x|^r)$. An exceptional curve $E_{r_j}$ corresponds to a vertex $\zeta(0, |x|^{r_j})$ which is mapped by $f_*$ into some annulus $U_j$. A closed point $E_{r_j} \cap E_{r_{j+1}}$ mapping to another closed point $E_{r_k} \cap E_{r_{k+1}}$ is observed as $f_*(U_j) \subseteq U_k$. An indeterminate point arises where $f_*(U_j)$ contains some $\zeta(0, |x|^{r_l})$; more specifically this means that $f(E_{r_j} \cap E_{r_{j+1}})$ contains $E_{r_l}$.

\subsection{Failed Stabilisation Through Graphs for Rational Maps}

Whilst under some conditions the method presented in \autoref{thm:blowupev} may work for (not just birational but) rational maps, it can \emph{easily} fail. Here I present a simple example continuing the work in this section.

\begin{prop}\label{prop:badgraphalg}
  Let $f : X \dashto X$ be the same map as above, on the Hirzebruch surface $X$. 
 If $\sigma_2 : (f_2, X_2) \to (f, X)$ is the point blowup of $(0, 2) \in X$, then $(f_2, X_2)$ is also algebraically stable. If however we apply $\pi_0 : (f_2, X_2) \to (g, Y)$, the blowdown of the proper transform of $\set{x=0}$, then $g$ is not AS and applying the smooth graph method of \autoref{thm:blowupev} will fail.
\end{prop}

\begin{proof}
In this section we already saw that $(f, X)$ is algebraically stable. The surface $X_2$ has one new exceptional curve, say $D$.  Just as we related $E_r$ to the family of curves $\gamma_r$ above, $D$ can be related to the family \[(s, t) \longmapsto (t, 2 + st).\]  Over this point $\sigma_1(D) = (2, 0)$, $D$ plays a similar role to the one $E_{\frac 11}$ did over the origin. Mapping this family forward by $f$ we get \[(s, t) \longmapsto \left(t^2, \frac{t^4 + (2 + st)^6}{(2 + st)^3}\right)\]
\[=(t^2, 8 + 12st + \bO(t^2)).\]
This new family sweeps out neither $E_0$ nor $D$ but targets the point $(0, 8) \in E_0$. The orbit of this point then continues along continuously defined points \[\left(0, 2^3\right) \mapsto \left(0, 2^{3^2}\right) \mapsto \left(0, 2^{3^3}\right) \mapsto \dots.\] Hence this is not destabilising and $f_2$ is algebraically stable.

Considering $(g, Y)$, it is clear that the proper transform of $D$ by $g$ is a point, namely the proper transform of $E_0$ by $\pi_0$; call this $P$. Now observe that \[g(P) = \pi_0 \circ f_2 \circ \pi_0^{-1}(P) = \pi_0 \circ f_2(E_0)\] (using that the composition is stable). Considering a similar calculation as above, one can see that $(0, 2^{\frac 13})$ is indeterminate with image $D$. Hence the total transform of $E_0$ by $f_2$ is $E_0 \cup D$. This shows that \[g(P) = \pi_0 \circ f_2(E_0) = D\] and hence $D \mapsto P \mapsto D$ is a destabilising orbit.

Suppose now we apply the smooth graph method to $(g, Y)$. Let $(g_1, \sgraph_g)$ be produced by the first iteration. Then $\alpha : \sgraph_g \to Y$ will blowup at least $P \in I(g)$, hence we can already factor $\alpha = \pi_0 \circ \alpha'$ with $\alpha' : \sgraph_g \to X_2$. We have seen earlier this section that $E_0, E_\frac 43 \mapsto E_0$; and just as $E_0$ has indeterminate points mapping to $D$, so does $E_\frac 43$. To see this, consider the following family of curves 
\[(s, t) \longmapsto (t^3, 2^{-\frac 13}t^4 + st^{10}).\]
The reader may verify that this is mapped by $f$ (or $g$) to a family of the form \[(s, t) \longmapsto (t^6, 2 - 3\cdot 2^{\frac43}st^6 + \bO(t^{12})),\] and this sweeps out $D$. Therefore for $\alpha$ to resolve this indeterminacy it must resolve the family of curves $(t^3, 2^{-\frac 13}t^4 + st^{10})$. To do this, it is necessary for $\alpha$ to blowup the origin (on $E_\frac 01$) to get $E_\frac 11$, then blowup further to produce $E_\frac 21$, $E_\frac 32$, $E_\frac 43$, and finally two more exceptional curves after this. In doing so we have created a surface which dominates the surface $X_0$ as in \autoref{thm:counterex}. Therefore by the theorem $(g_1, \sgraph_g)$ is not algebraically stable, moreover \emph{any} further blowups will result in an unstable rational map and thus the graph algorithm must fail.
\end{proof}
 
\subsection{Quadratic Example}

The following examples are based upon perhaps the simplest example of a rational but non-invertible map, $(x, y) \mapsto (x^2, y)$ on $\P^1 \times \P^1$. First we find a counter-example to the extension of \autoref{prop:factorblowup} to rational maps, as stated in \autoref{rmk:falseforrat}. Second, as mentioned in \autoref{rmk:untangledrat}, we give an example of an \nice{} rational map such that $\E(\alpha) \cap \E(\beta) \ne \emp$ in the smooth graph, hence making a generalised \autoref{prop:curveex} false.

For the following analysis, we consider a rational map $f : X \dashto Y$ with a possibly different domain and codomain, but each either equal to $\P^1 \times \P^1$ or satellite between $E_0$ and $E_1$. The irreducible curves over $\set{x=0}$ for these surfaces may be indexed separately by Farey parameters \mbox{$r_0 < \cdots < r_m$} and \mbox{$s_0 < \cdots < s_n$} respectively. We let $D$ denote the proper transform of $\set{y=0}$ on these surfaces. First note that \autoref{prop:mapexc} can easily be adapted to this situation by adjusting notation. Second, this result can be extended to understand the behaviour of $f$ with the closed points $E_{r_m} \cap D$ and $E_{s_n} \cap D$, which correspond precisely to the open intervals of Farey parameters $(r_m, \infty)$ and $(s_n, \infty)$ respectively. As in the previous examples, we need to determine how a rational map acts on a family of curves $\gamma_q : \P^1 \times \P^1 \dashto X$.

\begin{prop}
 Let $\phi : X \to \P^1 \times \P^1$, $\psi : Y \to \P^1 \times \P^1$ be either trivial maps or satellite between $E_0$ and $E_1$. Define the rational map $f : X \dashto Y$ to be the one induced by $(x, y) \mapsto (x^2, y)$ on $\P^1 \times \P^1$. Then $f \circ \gamma_{\frac ab} = \gamma_{\frac a{2b}}$ so the action of $f$ over $\set{x=0}$ is governed by $T_f(\frac ab) = \frac{a}{2b}$.
\end{prop}

\begin{ex}\label{ex:ratfactor}
 Let $Y = \P^1 \times \P^1$, $X$ be $\P^1 \times \P^1$ blown up at the origin, $(0, 0)$, and consider the lift $f : X \dashto Y$ of $(x, y) \mapsto (x^2, y)$. Then $f$ is a morphism, and $E_1 = \E(f)$ is contracted by $f$ to $p = (0, 0) \in I(f^{-1})$. However, if we blowup $p = (0, 0)$ as \autoref{prop:mapexc} suggests, with $\pi : Y' = X \to Y$ then $f' = \pi^{-1} \circ f$ still contracts $E_1 = \E(f')$ to a point, $E_0 \cap E_1$, and thus $\comp(f') = \comp(f) \ne \comp(f) - 1$. Worse, $f'$ is not even a morphism because $E_1\cap D \in I(f')$ is an indeterminate point whose image is $E_1$.
 
 Since $T_f(1) = \frac 12 \in (0, \infty)$, and $E_0$ is the only indexed curve on $Y$, by the extension of \autoref{prop:mapexc} discussed above, $E_1$ is mapped to a closed point corresponding to $E_0 \cap D = (0, 0)$. The only preimage of $0$ under $T_f$ is $0$, and since $E_0 \subset X$, there are no indeterminate points. After the blowup $\pi : Y' = X \to Y$ we have $T_{f'} = T_f$ and an extra curve $E_1 \subset Y$. We now find that $0 < T_f(1) = \frac 12 < 1$, meaning the image of $E_1$ is $E_0 \cap E_1$, and also $T_f((1, \infty)) \ni 1$ so the point $E_1\cap D$ maps to $E_1$. Note that in either case, the proper transform of $E_0$ under $f$ or $f'$ is $E_0$ and moreover $\E(f) = \E(f') = E_1$.
\end{ex}

\begin{ex}\label{ex:untangledrat}
 Let $X = \P^1 \times \P^1$, $Y$ be $\P^1 \times \P^1$ blown up at the origin, $(0, 0)$, and consider the lift $f : X \to Y$ of $(x, y) \mapsto (x^2, y)$. Then $(0, 0)$ is indeterminate but $f$ has no exceptional curves, hence $f$ is \nice{}. However the smooth graph, obtained by two blowups over $(0, 0)$, has a curve which is contracted by both $\alpha$ and $\beta$.

The challenge here is to blowup the domain $X$ just enough to resolve the indeterminacy of $f$, so that the resulting surface is the smooth graph of $f$. As in the previous example, $E_0 \subset Y$ is not generating indeterminacy since $T_f^{-1}(0) = \set 0$, but we do have $T_f^{-1}(1) = \set 2 \subset (0, \infty)$, hence $(0, 0)$ is indeterminate by \autoref{prop:mapexc}. Clearly to lift $f$ to a morphism $\beta : \sgraph_{\!f} \to Y$ through $\alpha : \sgraph_{\!f} \to X$, the surface $\sgraph_{\!f}$ must contain $E_2$, otherwise some closed point would map to $E_1$. Blowing up the origin gives the exceptional curve $E_1$, and one can check that blowing up $E_1 \cap D$ yields $E_2$. We claim this double blowup is the correct $\alpha : \sgraph_{\!f} \to X$, lifting $f$ to a morphism $\beta$, with $T_\beta = T_f$. Note that we have the curves $E_0, E_1, E_2$ on $\sgraph_{\!f}$, and $E_0, E_1$ in $Y$. Indeed, using \autoref{prop:mapexc} one can see there are no indeterminate points, since $T_\beta^{-1}(0) = \set 0$ and $T_\beta^{-1}(1) = \set 2$. Furthermore, $\beta(E_1) = E_0 \cap E_1$ because $T_\beta(1) = \frac 12 \in (0, 1)$, and therefore $E_1 \subseteq \E(\alpha)\cap\E(\beta) \ne \emp$.
\end{ex}
 
\bibliographystyle{alpha}
\bibliography{Stability}

\end{document}